\documentclass[11pt,twoside]{amsart}
\usepackage{xcolor}
\usepackage{ulem}
\usepackage{lipsum}
\usepackage{multirow}
\usepackage{pdflscape}
\usepackage{rotating}
\usepackage{graphicx}
\usepackage{amssymb}
\usepackage{amsmath}
\usepackage{amsthm}
\usepackage{amsfonts}
\usepackage{amsmath, amsfonts, amssymb}
\usepackage{enumerate}
\usepackage{enumitem}
\newtheorem{theorem}{Theorem}[section]

\newtheorem{example}[theorem]{Example}

\newtheorem{defn}[theorem]{Definition}
\numberwithin{equation}{section}

\begin{document}
\title{On certain three algebras generated by binary algebras}
\author{H. Ahmed$^{1,3}$, M.A.A. Ahmed$^{2,3}$, Sh.K. Said Husain$^{1,2}$, Witriany Basri$^{1}$}

\thanks{{\scriptsize
emails: houida\_m7@yahoo.com; mohammed\_h7@yahoo.com; kartini@upm.edu.my; witriany@upm.edu.my.}}
\maketitle
\begin{center}
\address{$^{1}$Department of Math., Faculty of Science, UPM, Selangor, Malaysia}\\
  \address{$^2$Institute for Mathematical Research (INSPEM), UPM, Serdang, Selangor, Malaysia}\\
\address{$^3$Depart. of Math., Faculty of Science, Taiz University, Taiz, Yemen}
\end{center}
\begin{abstract}   This paper's central theme is to prove the existence of an n-algebra whose multiplication cannot be expressed employing any binary operation. Furthermore, to prove if two algebras are not isomorphic, this property does not hold for $3$-algebras corresponding to these two algebras. The proof drives applying some results gotten early applying a new approach for the classification algebras problem, introduced recently, which showed great success in solving many classification algebras problems.
\end{abstract}
\pagestyle{myheadings}
\markboth{\rightline {\sl   H. Ahmed, et.al. }}
         {\leftline{\sl   On certain $3$-algebras generated by binary algebras}}

\bigskip
\newpage
\section{Introduction}
In 1969 \cite{6}, Kurosh introduced the notion of multilinear operator algebra. It is known that such algebraic structures are attractive for their applications to problems of modern mathematical physics. In 1973 \cite{7}, Nambu proposed an exciting generalization of classical Hamiltonian mechanics; the Nambu bracket is a generalization of the classical Poisson bracket.

Indeed, the advance of theoretical physics of quantum mechanics and the discovery of Nambu mechanics (see \cite{7}), together with Okubo's work on the Yang-Baxter equation (see \cite{8}), gave impetus to significant development on triple algebra ($3$-algebras).

Furthermore, Carlsson, Lister, and Loos have studied triple algebra of associative type (see \cite{51, 61, 62}). Hestenes provided the typical and founding example of totally associative triple algebra (see \cite{52}).

In this article, we give basic definitions and examples related to general $n$-algebras, and we shall focus our attention on $3$-algebras structures generated by binary algebras presented recently in \cite{1}. Then, we introduce the definition of totally associative $3$-algebras with examples, which show. Owing to the large size of the matrices involved in our computations of totally associative $3$-algebras, we present only Mathematica's results.

\section{Preliminaries}
     Let $\mathbb{F}$ be any field and the product $A\otimes B$ is the Kronecker product which stands for the matrix with blocks $(a_{ij}B),$ where $A=(a_{ij})$ and $B$ are matrices over $\mathbb{F}$.

\begin{defn} A vector space  $\mathbb{A}$ over $\mathbb{F}$ with multiplication $\cdot :\mathbb{A}\times \mathbb{A}\rightarrow \mathbb{A}$ given by $(\mathbf{u},\mathbf{v})\mapsto \mathbf{u}\cdot \mathbf{v}$ such that
\begin{itemize}
  \item $(\alpha\mathbf{u}+\beta\mathbf{v})\cdot \mathbf{w}=\alpha(\mathbf{u}\cdot \mathbf{w})+\beta(\mathbf{v}\cdot \mathbf{w}),$
  \item $\mathbf{w}\cdot (\alpha\mathbf{u}+\beta\mathbf{v})=\alpha(\mathbf{w}\cdot \mathbf{u})+\beta(\mathbf{w}\cdot \mathbf{v}),$
\end{itemize} whenever $\mathbf{u}, \mathbf{v}, \mathbf{w}\in \mathbb{A}$ and $\alpha, \beta\in \mathbb{F}$, is said to be an algebra.\end{defn}

\begin{defn} Two algebras $\mathbb{A}$ and $\mathbb{B}$ are called isomorphic if there is an invertible linear map  $\mathbf{f}:\mathbb{A}\rightarrow \mathbb{B} $ such that \begin{equation}\mathbf{f}(\mathbf{u}\cdot_{\mathbb{A}} \mathbf{v})=\mathbf{f}(\mathbf{u})\cdot_{\mathbb{B}} \mathbf{f}(\mathbf{v})\mbox{ whenever } \mathbf{u}, \mathbf{v}\in \mathbb{A}.\end{equation}\end{defn}
\begin{defn} A vector space $V$  over $\mathbb{F}$ equipped by a multilinear map  $f:\underbrace{V\times V\times ...\times V}_{n-times} \longrightarrow V$ is said to be a $n$-algebra, that means:
\begin{itemize}
  \item $f(x_1,x_2,..., x_i+x'_i,...,x_n)=f(x_1,x_2,..., x_i,...,x_n)+f(x_1,x_2,..., x'_i,...,x_n)$
  \item $f(x_1,x_2,..., \lambda x_i,...,x_n)= \lambda f(x_1,x_2,..., x_i,...,x_n)$
\end{itemize} where $(x_1,x_2,..., x_i,...,x_n) \in V.$\end{defn}
\begin{example}\label{1} Let $\mathbb{A}=(V,\mu)$ be an algebra over a field $\mathbb{F}$: Then multilinear map
\begin{equation}\label{eq1}f(x_1,x_2,…,x_n)=\underbrace{\mu(x_1,\mu(x_2,..., \mu(x_{n-1},x_n)...))}_{(n-1)-times}\end{equation}
defines an $n$-algebra structure on $V$.\end{example}

\section{Classification approach of $m$-dimensional $3$-algebras}
  Let $\mathbb{A}$ be $m$-dimensional $3$-algebra over $\mathbb{F}$ and $\mathbf{e}=(e^1,e^2,...,e^m)$ its basis. Then the multilinear map $\cdot$ is represented by a matrix $A=(A^l_{ijk})\in M(m\times m^3;\mathbb{F})$ as follows
   \begin{equation}\label{bi} \mathbf{u}\cdot \mathbf{v} \cdot \mathbf{w} =\mathbf{e}A(u\otimes v \otimes w),\end{equation}
    for $\mathbf{u}=\mathbf{e}u,\mathbf{v}=\mathbf{e}v,\mathbf{w}=\mathbf{e}w,$ where $u = (u_1, u_2, ..., u_m)^T,$  $v = (v_1, v_2, ..., v_m)^T$ and $w = (w_1, w_2, ..., w_m)^T$ are column coordinate vectors of $\mathbf{u},$ $\mathbf{v},$ and $\mathbf{w},$ respectively.
The matrix $A\in M(m\times m^3;\mathbb{F})$ defined above is called the matrix of structural constants (MSC) of $\mathbb{A}$ with respect to the basis $\mathbf{e}$. Further we assume that a basis $\mathbf{e}$ is fixed and we do not make a difference between the algebra
$\mathbb{A}$ and its MSC $A$.

If $\mathbf{e}'=(e'^1,e'^2,...,e'^m)$ is another basis of $\mathbb{A}$, $\mathbf{e}'g=\mathbf{e}$ with $g\in G=GL(m;\mathbb{F})$, and  $A'$ is MSC of $\mathbb{A}$ with respect to $\mathbf{e}'$ then it is known that
 \begin{equation} \label{2} A'=gA(g^{-1})^{\otimes 3}\end{equation}
 is valid (see \cite{5}). Thus, we can reformulate the isomorphism of algebras as follows.
  \begin{defn} Two $m$-dimensional $3$-algebras $\mathbb{A}$, $\mathbb{B}$ over $\mathbb{F}$, given by
	their matrices of structure constants $A$, $B$, are said to be isomorphic if $B=gA(g^{-1})^{\otimes 3}$ holds true for some $g\in GL(m;\mathbb{F})$.\end{defn}

Further we consider only the case $m=2$ then $\mathbb{A}$ can be represented by its matrix of structural constants (MSC) $A=(\gamma_{ijk}^{l} )\in M(2\times 2^3;F)$ where $i,j,k ,l=1,2 $ and $\gamma_{ijk}^{l} \in \mathbb{F},$ as follows:
\[A=\left(
\begin{array}{cccccccc}
 \gamma_{111}^{1} & \gamma_{112}^{1}& \gamma_{121}^{1} & \gamma_{122}^{1} & \gamma_{211}^{1} & \gamma_{212}^{1} & \gamma_{221}^{1} & \gamma_{222}^{1} \\
 \gamma_{111}^{2} & \gamma_{112}^{2}& \gamma_{121}^{2} & \gamma_{122}^{2} & \gamma_{211}^{2} & \gamma_{212}^{2} & \gamma_{221}^{2} & \gamma_{222}^{2} \\
\end{array}
\right)\]

(for more information refer to \cite{2}).
\section{$3$-algebras generated by binary algebras}
Due to \cite{1} we have the following classification theorems according to $Char(\mathbb{F})\neq 2,3.$

  \begin{theorem}\label{theorem1} Over an algebraically closed field $\mathbb{F}$ $(Char(\mathbb{F})\neq 2$ and $3)$, any non-trivial $2$-dimensional algebra is isomorphic to only one of the following algebras listed by their matrices of structure constants:
	\begin{itemize}
	\item $A_{1}(\mathbf{c})=\left(
	\begin{array}{cccc}
	\alpha_1 & \alpha_2 &\alpha_2+1 & \alpha_4 \\
	\beta_1 & -\alpha_1 & -\alpha_1+1 & -\alpha_2
	\end{array}\right),\ \mbox{where}\ \mathbf{c}=(\alpha_1, \alpha_2, \alpha_4, \beta_1)\in \mathbb{F}^4,$
	\item $A_{2}(\mathbf{c})=\left(
	\begin{array}{cccc}
	\alpha_1 & 0 & 0 & 1 \\
	\beta _1& \beta _2& 1-\alpha_1&0
	\end{array}\right)\simeq \left(
	\begin{array}{cccc}
	\alpha_1 & 0 & 0 & 1 \\
	-\beta _1& \beta _2& 1-\alpha_1&0
	\end{array}\right),\ \mbox{where}\ \mathbf{c}=(\alpha_1, \beta_1, \beta_2)\in \mathbb{F}^3,$
	\item $A_{3}(\mathbf{c})=\left(
	\begin{array}{cccc}
	0 & 1 & 1 & 0 \\
	\beta _1& \beta _2 & 1&-1
	\end{array}\right),\ \mbox{where}\ \mathbf{c}=(\beta_1, \beta_2)\in \mathbb{F}^2,$
	\item $A_{4}(\mathbf{c})=\left(
	\begin{array}{cccc}
	\alpha _1 & 0 & 0 & 0 \\
	0 & \beta _2& 1-\alpha _1&0
	\end{array}\right),\ \mbox{where}\ \mathbf{c}=(\alpha_1, \beta_2)\in \mathbb{F}^2,$
	\item $A_{5}(\mathbf{c})=\left(
	\begin{array}{cccc}
	\alpha_1& 0 & 0 & 0 \\
	1 & 2\alpha_1-1 & 1-\alpha_1&0
	\end{array}\right),\ \mbox{where}\ \mathbf{c}=\alpha_1\in \mathbb{F},$
	\item $A_{6}(\mathbf{c})=\left(
	\begin{array}{cccc}
	\alpha_1 & 0 & 0 & 1 \\
	\beta _1& 1-\alpha_1 & -\alpha_1&0
	\end{array}\right)\simeq \left(
	\begin{array}{cccc}
	\alpha_1 & 0 & 0 & 1 \\
	-\beta _1& 1-\alpha_1 & -\alpha_1&0
	\end{array}\right),\ \mbox{where}\ \mathbf{c}=(\alpha_1, \beta_1)\in \mathbb{F}^2,$
	\item $A_{7}(\mathbf{c})=\left(
	\begin{array}{cccc}
	0 & 1 & 1 & 0 \\
	\beta_1& 1& 0&-1
	\end{array}\right),\ \mbox{where}\ \mathbf{c}=\beta_1\in \mathbb{F},$
	\item $A_{8}(\mathbf{c})=\left(
	\begin{array}{cccc}
	\alpha_1 & 0 & 0 & 0 \\
	0 & 1-\alpha_1 & -\alpha_1&0
	\end{array}\right),\ \mbox{where}\ \mathbf{c}=\alpha_1\in \mathbb{F},$
\item $A_{9}=\left(
	\begin{array}{cccc}
	\frac{1}{3}& 0 & 0 & 0 \\
	1 & \frac{2}{3} & -\frac{1}{3}&0
	\end{array}\right),$
	\item $A_{10}=\left(
	\begin{array}{cccc}
	0 & 1 & 1 & 0 \\
	0 &0& 0 &-1
	\end{array}
	\right),$
	\item $A_{11}=\left(
	\begin{array}{cccc}
	0 & 1 & 1 & 0 \\
	1 &0& 0 &-1
	\end{array}
	\right),$
	\item $A_{12}=\left(
	\begin{array}{cccc}
	0 & 0 & 0 & 0 \\
	1 &0&0 &0\end{array}
	\right).$
\end{itemize}\end{theorem}
\begin{example}\label{2} In example \ref{1}, if $V $ is $2$-dimensional vector space with a fixed basis $\{e_1,e_2\}$ and $(V,f)$ is $3$-algebra structure on $V.$ Then $f$ and $\mu$ can be expressed by their structure constants as follows: $f(e_i,e_j,e_k )=\gamma_{ijk}^{1}e_1+\gamma_{ijk}^{2} e_2$ and $\mu (e_r,e_s )=\eta_{rs}^1 e_1+\eta_{rs}^2 e_2.$ Then due to (\ref{eq1}) we get the system of equations
\begin{equation}\label{eq2}\left.\begin{array}{cc}
                             \gamma_{ijk}^{1}= & \eta_{jk}^1\eta_{i1}^1+\eta_{jk}^2\eta_{i2}^1 \\
                             \gamma_{ijk}^{2}= & \eta_{jk}^1\eta_{i1}^2+\eta_{jk}^2\eta_{i2}^2
                           \end{array}\right\}
\end{equation}
we get $2^4 = 16$ equations for $2^3 = 8$
unknowns (the coefficients $\eta_{ij}^k$), which cannot be solved in general, except maybe for some very special cases.\\
Indeed, using (\ref{eq2}), we can find the $3$-algebras corresponding to all algebras presented in \cite{1} under this procedure (see Table \ref{T3algebras}). \end{example}

%\begin{adjustbox}{angle=90}
%\rotatebox{90}{
\begin{tiny}
\begin{sidewaystable}
\vspace*{15cm}
\begin{tabular}{|c|c|}
  \hline
  algebra & $3$-algebra \\
  \hline
  $A_1(\mathbf{c})%=\left(
%\begin{array}{cccc}
% \alpha _1 & \alpha _2 & \alpha _2+1 & \alpha _4 \\
 %\beta _1 & -\alpha _1 & 1-\alpha _1 & -\alpha _2 \\
%\end{array}
%\right)
$& $B_1(\mathbf{c})=\left(
\begin{array}{cccccccc}
 \alpha _2 \beta _1+\alpha _1^2 & 0 & \alpha _1+\alpha _2 & \alpha _1 \alpha _4-\alpha _2^2 & \alpha _4 \beta _1+\alpha _2 \alpha _1+\alpha _1 & \alpha _2^2-\alpha _2-\alpha _1 \alpha _4 & \alpha _2^2+2 \alpha _2-\alpha _1 \alpha _4+\alpha _4+1 & \alpha _4 \\
 0 & \alpha _2 \beta _1+\alpha _1^2 & \alpha _2 \beta _1+\alpha _1^2-\alpha _1+\beta _1 & \alpha _4 \beta _1+\alpha _1 \alpha _2 & -\alpha _2 \beta _1-\alpha _1^2+\alpha _1 & \alpha _2 & 1-\alpha _1 & \alpha _2^2-\alpha _1 \alpha _4+\alpha _4 \\
\end{array}
\right)$ \\
  \hline
   $A_2(\alpha _1,\beta _1,\beta _2)%=\left(
%\begin{array}{cccc}
% \alpha _1 & 0 & 0 & 1 \\
% \beta _1 & \beta _2 & 1-\alpha _1 & 0 \\
%\end{array}
%\right)
$ & $B_2(\alpha _1,\beta _1,\beta _2)=\left(
\begin{array}{cccccccc}
 \alpha _1^2 & 0 & 0 & \alpha _1 & \beta _1 & \beta _2 & 1-\alpha _1 & 0 \\
 \alpha _1 \beta _1+\beta _2 \beta _1 & \beta _2^2 & \left(1-\alpha _1\right) \beta _2 & \beta _1 & \alpha _1 \left(1-\alpha _1\right) & 0 & 0 & 1-\alpha _1 \\
\end{array}
\right)$ \\
  \hline
  $A_3(\beta _1,\beta _2)%=\left(
%\begin{array}{cccc}
 %0 & 1 & 1 & 0 \\
 %\beta _1 & \beta _2 & 1 & -1 \\
%\end{array}
%\right)
$ & $B_3(\beta _1,\beta _2)=\left(
\begin{array}{cccccccc}
 \beta _1 & \beta _2 & 1 & -1 & 0 & 1 & 1 & 0 \\
 \beta _1 \beta _2 & \beta _2^2+\beta _1 & \beta _1+\beta _2 & -\beta _2 & -\beta _1 & 1-\beta _2 & 0 & 1 \\
\end{array}
\right)$ \\
  \hline
  $A_4(\alpha _1,\beta _2)%=\left(
%\begin{array}{cccc}
% \alpha _1 & 0 & 0 & 0 \\
% 0 & \beta _2 & 1-\alpha _1 & 0 \\
%\end{array}
%\right)
$ & $B_4(\alpha _1,\beta _2)=\left(
\begin{array}{cccccccc}
 \alpha _1^2 & 0 & 0 & 0 & 0 & 0 & 0 & 0 \\
 0 & \beta _2^2 & \left(1-\alpha _1\right) \beta _2 & 0 & \alpha _1 \left(1-\alpha _1\right) & 0 & 0 & 0 \\
\end{array}
\right) $\\
  \hline
  $A_5(\alpha _1)%=\left(
%\begin{array}{cccc}
 %\alpha _1 & 0 & 0 & 0 \\
 %1 & 2 \alpha _1-1 & 1-\alpha _1 & 0 \\
%\end{array}
%\right)
$ & $B_5(\alpha _1)=\left(
\begin{array}{cccccccc}
 \alpha _1^2 & 0 & 0 & 0 & 0 & 0 & 0 & 0 \\
 3 \alpha _1-1 & \left(2 \alpha _1-1\right){}^2 & \left(2 \alpha _1-1\right) \left(1-\alpha _1\right) & 0 & \alpha _1 \left(1-\alpha _1\right) & 0 & 0 & 0 \\
\end{array}
\right) $\\
  \hline
   $A_6(\alpha _1,\beta _1)%=\left(
%\begin{array}{cccc}
% \alpha _1 & 0 & 0 & 1 \\
 %\beta _1 & 1-\alpha _1 & -\alpha _1 & 0 \\
%\end{array}
%\right)
$ &$ B_6(\alpha _1,\beta _1)=\left(
\begin{array}{cccccccc}
 \alpha _1^2 & 0 & 0 & \alpha _1 & \beta _1 & 1-\alpha _1 & -\alpha _1 & 0 \\
 \beta _1 & \left(1-\alpha _1\right){}^2 & -\alpha _1 \left(1-\alpha _1\right) & \beta _1 & -\alpha _1^2 & 0 & 0 & -\alpha _1 \\
\end{array}
\right)$ \\
  \hline
  $A_7(\beta _1)%=\left(
%\begin{array}{cccc}
% 0 & 1 & 1 & 0 \\
% \beta _1 & 1 & 0 & -1 \\
%\end{array}
%\right)
$ & $B_7(\beta _1)=\left(
\begin{array}{cccccccc}
 \beta _1 & \beta _1+1 & 0 & -1 & 0 & 1 & 1 & 0 \\
 \beta _1 & 1 & \beta _1 & -1 & -\beta _1 & -1 & 0 & 1 \\
\end{array}
\right)$ \\
  \hline
   $A_8(\alpha _1)%=\left(
%\begin{array}{cccc}
% \alpha _1 & 0 & 0 & 0 \\
% 0 & 1-\alpha _1 & -\alpha _1 & 0 \\
%\end{array}
%\right)
$ & $B_8(\alpha _1)=\left(
\begin{array}{cccccccc}
 \alpha _1^2 & 0 & 0 & 0 & \alpha _1^2 & 0 & 0 & 0 \\
 0 & \left(1-\alpha _1\right){}^2 & -\alpha _1 \left(1-\alpha _1\right) & 0 & 0 & 0 & 0 & 0 \\
\end{array}
\right)$ \\
  \hline
  $A_9%=\left(
%\begin{array}{cccc}
 %\frac{1}{3} & 0 & 0 & 0 \\
% 1 & \frac{2}{3} & -\frac{1}{3} & 0 \\
%\end{array}
%\right)
$ & $B_9=\left(
\begin{array}{cccccccc}
 \frac{1}{9} & 0 & 0 & 0 & 0 & 0 & 0 & 0 \\
 1 & \frac{4}{9} & -\frac{2}{9} & 0 & -\frac{1}{9} & 0 & 0 & 0 \\
\end{array}
\right)$ \\
  \hline
 $ A_{10}%=\left(
%\begin{array}{cccc}
% 0 & 1 & 1 & 0 \\
% 0 & 0 & 0 & -1 \\
%\end{array}
%\right)
$ & $B_{10}=\left(
\begin{array}{cccccccc}
 0 & 0 & 0 & -1 & 0 & 1 & 1 & 0 \\
 0 & 0 & 0 & 0 & 0 & 0 & 0 & 1 \\
\end{array}
\right)$ \\
  \hline
 $A_{11}%=\left(
%\begin{array}{cccc}
% 0 & 1 & 1 & 0 \\
% 1 & 0 & 0 & -1 \\
%\end{array}
%\right)
$ & $B_{11}=\left(
\begin{array}{cccccccc}
 1 & 0 & 0 & -1 & -1 & 1 & 1 & 0 \\
 0 & 1 & 1 & 0 & 0 & 0 & 0 & 1 \\
\end{array}
\right)$ \\
  \hline
 $A_{12}%=\left(
%\begin{array}{cccc}
% 0 & 0 & 0 & 0 \\
% 1 & 0 & 0 & 0 \\
%\end{array}
%\right)
$ & trivial \\
  \hline
\end{tabular}\caption{$2$-dimensional $3$-algebras generated by binary algebras(where $\mathbf{c}=(\alpha_1, \alpha_2, \alpha_4, \beta_1)$)}
\label{T3algebras}\end{sidewaystable}\end{tiny}
From the table we can see
 \begin{itemize}
 \item $\left(
\begin{array}{cccccccc}
 1 & 0 & 0 & 1 & 0 & 1 & -1 & 0 \\
 0 & -1 & 1 & 0 & 1 & 0 & 0 & 1 \\
\end{array}
\right)$ is a $3$-algebra which cannot be expressed by any algebras in the above theorem and it is not isomorphic to any $B_i $ where $i=1,...,11.$
\item On the other hand, $A_4 (1/3,-1/3) $ and $A_5 (1/3)$ are not isomorphic algebras, but from these two non-isomorphic algebras, we get one $3$-algebra $\left(
\begin{array}{cccccccc}
 \frac{1}{9} & 0 & 0 & 0 & 0 & 0 & 0 & 0 \\
 0 & \frac{1}{9} & -\frac{2}{9} & 0 & \frac{2}{9} & 0 & 0 & 0 \\
\end{array}
\right).$
\item Also $A_4 (1,-1) $ and $A_4(1,1)$ are not isomorphic algebras but from these two algebras we get one $3$-algebra $\left(
\begin{array}{cccccccc}
 1 & 0 & 0 & 0 & 0 & 0 & 0 & 0 \\
 0 & 1 & 0 & 0 & 0 & 0 & 0 & 0 \\
\end{array}
\right).$\end{itemize}
\section{Totally associative $3$-algebras}
In this section, we introduce the following definition of totally associative $3$-algebra using its MSC.
\begin{defn}\label{3ass}
 A $3$-algebra $\mathbb{A}$ is a totally associative $3$-algebra if
 \begin{equation}\label{3oas}
   (\mathbf{u}\cdot \mathbf{v}\cdot\mathbf{w})\cdot \mathbf{x}\cdot \mathbf{y}=
   \mathbf{u} \cdot(\mathbf{v}\cdot \mathbf{w}\cdot \mathbf{x})\cdot \mathbf{y}=
   \mathbf{u} \cdot\mathbf{v}\cdot (\mathbf{w}\cdot \mathbf{x}\cdot \mathbf{y})
 \end{equation} for all $\mathbf{u}, \mathbf{v}, \mathbf{w}, \mathbf{x}, \mathbf{y} \in \mathbb{A}.$
\end{defn}
\begin{example}
  Let $\{e_1, e_2\}$ be a basis of a $2$-dimensional $3$-algebra $\mathbb{A},$ the multilinear map "$\cdot$" given by:
$$ e_1 \cdot e_1\cdot e_1=e_1, \ \ e_1\cdot e_1\cdot e_2=e_2$$ defines a totally associative $3$-algebra.
\end{example}
According to (\ref{3oas}) and (\ref{bi}), we can reformulate Definition \ref{3ass} as follows.
\begin{defn}\label{n3ass}
  A $2$-dimensional $3$-algebra $\mathbb{A}$ with multiplication $\cdot$ over a field $\mathbb{F}$ is said to be a totally associative $3$-algebra if all of the following conditions are met:
  \begin{subequations}\label{n3asso}
\begin{align}
A\left(A\otimes I\otimes I-I\otimes A \otimes I\right)=0,\label{n3asso:subeq_a}\\
 A\left(A\otimes I\otimes I-I\otimes I \otimes A\right)=0,\label{n3asso:subeq_b}\\
A\left(I \otimes A \otimes I-I\otimes I \otimes A \right)=0,\label{n3asso:subeq_c}
  \end{align}\end{subequations}

where $I$ is the identity $2\times 2$ matrix.\end{defn}

Using a computation program (here, we use Mathematica), it is easy to verify that the $3$-algebras from the list in Table (\ref{T3algebras}) satisfying the system (\ref{n3asso}) are:
\begin{enumerate}[label=(\roman*)]
\item $B_2(\alpha_1,\beta _1,\beta _2)$ when
\begin{itemize}
\item $\alpha _1=0, \ \beta _1=0, \ \beta _2=0,$
\item $\alpha _1=\frac{1}{2}, \ \beta _1=0, \ \beta _2=-\frac{1}{2},$
\item $\alpha _1=\frac{1}{2}, \  \beta _1=0, \ \beta _2=\frac{1}{2}.$\end{itemize}

\item $B_4(\alpha_1,\beta _2)$ when
\begin{itemize}
\item $\alpha _1=0, \ \beta _2=0,$
\item $\alpha _1=\frac{1}{2}, \ \beta _2=-\frac{1}{2},$
\item $\alpha _1=\frac{1}{2}, \ \beta _2=0,$
\item $\alpha _1=\frac{1}{2}, \ \beta _2=\frac{1}{2},$
\item $\alpha _1=1, \ \beta _2=-1,$
\item $\alpha _1=1, \ \beta _2=0,$
\item $\alpha _1=1, \ \beta _2=1.$\end{itemize}\end{enumerate}
That means from the above list we get the following totally associative $3$-algebras:
\begin{enumerate}[label=(\roman*)]
\item $B_2(0,0,0)=\left(
\begin{array}{cccccccc}
 0 & 0 & 0 & 0 & 0 & 0 & 1 & 0 \\
 0 & 0 & 0 & 0 & 0 & 0 & 0 & 1 \\
\end{array}
\right) ,$
\item $ B_2(\frac{1}{2},0,-\frac{1}{2})=\left(
\begin{array}{cccccccc}
 \frac{1}{4} & 0 & 0 & \frac{1}{2} & 0 & -\frac{1}{2} & \frac{1}{2} & 0 \\
 0 & \frac{1}{4} & \frac{1}{4} & 0 & \frac{1}{4} & 0 & 0 & \frac{1}{2} \\
\end{array}
\right),$
\item $ B_2(\frac{1}{2},0,\frac{1}{2})=\left(
\begin{array}{cccccccc}
 \frac{1}{4} & 0 & 0 & \frac{1}{2} & 0 & \frac{1}{2} & \frac{1}{2} & 0 \\
 0 & \frac{1}{4} & \frac{1}{4} & 0 & \frac{1}{4} & 0 & 0 & \frac{1}{2} \\
\end{array}
\right)$
\item $B_4(\frac{1}{2},-\frac{1}{2})=\left(
\begin{array}{cccccccc}
 \frac{1}{4} & 0 & 0 & 0 & 0 & 0 & 0 & 0 \\
 0 & \frac{1}{4} & -\frac{1}{4} & 0 & \frac{1}{4} & 0 & 0 & 0 \\
\end{array}
\right) ,$
\item $B_4(\frac{1}{2},0)=\left(
\begin{array}{cccccccc}
 \frac{1}{4} & 0 & 0 & 0 & 0 & 0 & 0 & 0 \\
 0 & 0 & 0 & 0 & \frac{1}{4} & 0 & 0 & 0 \\
\end{array}
\right),$
\item $ B_4(\frac{1}{2},\frac{1}{2})=\left(
\begin{array}{cccccccc}
 \frac{1}{4} & 0 & 0 & 0 & 0 & 0 & 0 & 0 \\
 0 & \frac{1}{4} & \frac{1}{4} & 0 & \frac{1}{4} & 0 & 0 & 0 \\
\end{array}
\right),$
\item $B_4(1,-1)=B_4(1,1)= \left(
\begin{array}{cccccccc}
 1 & 0 & 0 & 0 & 0 & 0 & 0 & 0 \\
 0 & 1 & 0 & 0 & 0 & 0 & 0 & 0 \\
\end{array}
\right),$
\item  $B_4(1,0)= \left(
\begin{array}{cccccccc}
 1 & 0 & 0 & 0 & 0 & 0 & 0 & 0 \\
 0 & 0 & 0 & 0 & 0 & 0 & 0 & 0 \\
\end{array}
\right).$\end{enumerate}

According to (\cite{41}),  over an algebraically closed field, $\mathbb{F}$ of characteristic, not $2,\ 3$ every nontrivial $2$-dimensional associative algebra is isomorphic to only one algebra listed below by
\begin{enumerate}[label=(\roman*)]
\item $A_{2}\left(\displaystyle\frac{1}{2},0,\displaystyle\frac{1}{2}\right)=\left(
\begin{array}{cccc}
 \displaystyle\frac{1}{2} &0& 0 & 1 \\
 0 & \displaystyle\frac{1}{2} & \displaystyle\frac{1}{2} & 0
\end{array}
\right),$
\item $ A_{4}\left(1,0\right)=\left(
\begin{array}{cccc}
 1 & 0 & 0 & 0 \\
 0 & 0 & 0 & 0
\end{array}
\right),$
\item $ A_{4}\left(\displaystyle\frac{1}{2},\displaystyle\frac{1}{2}\right)=\left(
\begin{array}{cccc}
 \displaystyle\frac{1}{2} & 0 & 0 & 0 \\
 0 & \displaystyle\frac{1}{2} & \displaystyle\frac{1}{2} & 0
\end{array}
\right),$
\item $ A_{4}\left(1,1\right)=\left(
\begin{array}{cccc}
 1 & 0 & 0 & 0 \\
 0 & 1 & 0 & 0
\end{array}
\right),$
\item $ A_{4}\left(\displaystyle\frac{1}{2},0\right)=\left(
\begin{array}{cccc}
 \displaystyle\frac{1}{2} & 0 & 0 & 0 \\
 0 & 0 & \displaystyle\frac{1}{2} & 0
\end{array}
\right),$
\item $ A_{12}=\left(
\begin{array}{cccc}
 0 & 0 & 0 & 0 \\
 1 & 0 & 0 & 0
\end{array}
\right).$
\end{enumerate}
We conclude that there are totally associative $3$-algebras, $B_2(0,0,0),$ $ B_2(\frac{1}{2},0,-\frac{1}{2}),$ $B_4(\frac{1}{2},-\frac{1}{2}),$ $B_4(1,-1),$ are generated by non-associative algebras, $A_2(0,0,0),$ $ A_2(\frac{1}{2},0,-\frac{1}{2}),$ $A_4(\frac{1}{2},-\frac{1}{2}),$ $A_4(1,-1),$ respectively.
 \begin{center}{\textbf{Conclusion}}
\end{center}
Depending on the approach introduced in \cite{5} and applied in \cite{1}, one can study the classification of $n$-algebras and then study some identities of $n$-algebras (refer to \cite{2, 3, 4, 41,42, 43}).

\begin{center}{\textbf{Acknowledgement}}
\end{center}

This research was funded by Grant 01-01-18-2032FR Kementerian Pendidikan Malaysia.

\end{document}